\documentstyle[11pt]{article}
\begin{document}
%
%
%
\newtheorem{theorem}      {Th\'eor\`eme}[section]
\newtheorem{theorem*}     {theorem}
\newtheorem{proposition}  [theorem]{Proposition}
\newtheorem{definition}   [theorem]{Definition}
\newtheorem{e-lemme}        [theorem]{Lemma}
\newtheorem{cor}   [theorem]{Corollaire}
\newtheorem{resultat}     [theorem]{R\'esultat}
\newtheorem{eexercice}    [theorem]{Exercice}
\newtheorem{rrem}    [theorem]{Remarque}
\newtheorem{pprobleme}    [theorem]{Probl\`eme}
\newtheorem{eexemple}     [theorem]{Exemple}
\newcommand{\preuve}      {\paragraph{Preuve}}
\newenvironment{probleme} {\begin{pprobleme}\rm}{\end{pprobleme}}
\newenvironment{remarque} {\begin{rremarque}\rm}{\end{rremarque}}
\newenvironment{exercice} {\begin{eexercice}\rm}{\end{eexercice}}
\newenvironment{exemple}  {\begin{eexemple}\rm}{\end{eexemple}}
%
%
\newtheorem{e-theo}      [theorem]{Theorem}
\newtheorem{theo*}     [theorem]{Theorem}
\newtheorem{e-pro}  [theorem]{Proposition}
\newtheorem{e-def}   [theorem]{Definition}
\newtheorem{e-lem}        [theorem]{Lemma}
\newtheorem{e-cor}   [theorem]{Corollary}
\newtheorem{e-resultat}     [theorem]{Result}
\newtheorem{ex}    [theorem]{Exercise}
\newtheorem{e-rem}    [theorem]{Remark}
\newtheorem{prob}    [theorem]{Problem}
\newtheorem{example}     [theorem]{Example}
\newcommand{\proof}         {\paragraph{Proof~: }}
\newcommand{\hint}          {\paragraph{Hint}}
\newcommand{\heuristicproof}{\paragraph{heuristic proof}}
\newenvironment{e-probleme} {\begin{e-pprobleme}\rm}{\end{e-pprobleme}}
\newenvironment{e-remarque} {\begin{e-rremarque}\rm}{\end{e-rremarque}}
\newenvironment{e-exercice} {\begin{e-eexercice}\rm}{\end{e-eexercice}}
\newenvironment{e-exemple}  {\begin{e-eexemple}\rm}{\end{e-eexemple}}
\def \Re       {{\rm Re} }
\def \Im       {{\rm Im} }
\newcommand{\reell}    {{{\rm I\! R}^l}}
\newcommand{\reeln}    {{{\rm I\! R}^n}}
\newcommand{\reelk}    {{{\rm I\! R}^k}}
\newcommand{\reelm}    {{{\rm I\! R}^m}}
\newcommand{\reelp}    {{{\rm I\! R}^p}}
\newcommand{\reeld}    {{{\rm I\! R}^d}}
\newcommand{\reeldd}   {{{\rm I\! R}^{d\times d}}}
\newcommand{\reelnn}   {{{\rm I\! R}^{n\times n}}}
\newcommand{\reelnd}   {{{\rm I\! R}^{n\times d}}}
\newcommand{\reeldn}   {{{\rm I\! R}^{d\times n}}}
\newcommand{\reelkd}   {{{\rm I\! R}^{k\times d}}}
\newcommand{\reelkl}   {{{\rm I\! R}^{k\times l}}}
\newcommand{\reelN}    {{{\rm I\! R}^N}}
\newcommand{\reelM}    {{{\rm I\! R}^M}}
\newcommand{\reelplus} {{{\rm I\! R}^+}}
\newcommand{\reelo}    {{{\rm I\! R}\setminus\{0\}}}
\newcommand{\reld}    {{{\rm I\! R}_d}}
\newcommand{\relplus} {{{\rm I\! R}_+}}
\newcommand{\1}        {{\bf 1}}

\newcommand{\cov}      {{\hbox{cov}}}
\newcommand{\sss}      {{\cal S}}
\newcommand{\indic}    {{{\rm I\!\! I}}}
\newcommand{\pp}       {{{\rm I\!\!\! P}}}
\newcommand{\qq}       {{{\rm I\!\!\! Q}}}
\newcommand{\ee}       {{{\rm I\! E}}}

\newcommand{\B}        {{{\rm I\! B}}}
\newcommand{\cc}       {{{\rm I\!\!\! C}}}
\renewcommand{\pp}        {{{\rm I\!\!\! P}}}
\newcommand{\HHH}        {{{\rm I\! H}}}
\newcommand{\N}        {{{\rm I\! N}}}
\newcommand{\R}        {{{\rm I\! R}}}
\newcommand{\D}        {{{\rm I\! D}}}
\newcommand{\Z}       {{{\rm Z\!\! Z}}}
\newcommand{\C}        {{\bf C}}        
\newcommand{\T}        {{\bf T}}        
\newcommand{\E}        {{\bf E}}        
\newcommand{\rfr}[1]    {\stackrel{\circ}{#1}}
\newcommand{\equiva}    {\displaystyle\mathop{\simeq}}
\newcommand{\eqdef}     {\stackrel{\triangle}{=}}
\newcommand{\limps}     {\mathop{\hbox{\rm lim--p.s.}}}
\newcommand{\Limsup}    {\mathop{\overline{\rm lim}}}
\newcommand{\Liminf}    {\mathop{\underline{\rm lim}}}
\newcommand{\Inf}       {\mathop{\rm Inf}}
\newcommand{\vers}      {\mathop{\;{\rightarrow}\;}}
\newcommand{\versup}    {\mathop{\;{\nearrow}\;}}
\newcommand{\versdown}  {\mathop{\;{\searrow}\;}}
\newcommand{\vvers}     {\mathop{\;{\longrightarrow}\;}}
\newcommand{\cvetroite} {\mathop{\;{\Longrightarrow}\;}}
\newcommand{\ieme}      {\hbox{i}^{\hbox{\smalltype\`eme}}}
\newcommand{\eqps}      {\, \buildrel \rm \hbox{\rm\smalltype p.s.} \over \,}
\newcommand{\eqas}      {\,\buildrel\rm\hbox{\rm\smalltype a.s.} \over = \,}
\newcommand{\argmax}    {\hbox{{\rm Arg}}\max}
\newcommand{\argmin}    {\hbox{{\rm Arg}}\min}
\newcommand{\indep}{\perp\!\!\!\!\perp}
\newcommand{\abs}[1]{\left| #1 \right|}
\newcommand{\crochet}[2]{\langle #1 \,,\, #2 \rangle}
\newcommand{\espc}[3]   {E_{#1}\left(\left. #2 \right| #3 \right)}
\newcommand{\rang}{\hbox{rang}}
\newcommand{\rank}{\hbox{rank}}
\newcommand{\signe}{\hbox{signe}}
\newcommand{\sign}{\hbox{sign}}

\newcommand\hA{{\widehat A}}
\newcommand\hB{{\widehat B}}
\newcommand\hC{{\widehat C}}
\newcommand\hD{{\widehat D}}
\newcommand\hE{{\widehat E}}
\newcommand\hF{{\widehat F}}
\newcommand\hG{{\widehat G}}
\newcommand\hH{{\widehat H}}
\newcommand\hI{{\widehat I}}
\newcommand\hJ{{\widehat J}}
\newcommand\hK{{\widehat K}}
\newcommand\hL{{\widehat L}}
\newcommand\hM{{\widehat M}}
\newcommand\hN{{\widehat N}}
\newcommand\hO{{\widehat O}}
\newcommand\hP{{\widehat P}}
\newcommand\hQ{{\widehat Q}}
\newcommand\hR{{\widehat R}}
\newcommand\hS{{\widehat S}}
\newcommand\hTT{{\widehat T}}
\newcommand\hU{{\widehat U}}
\newcommand\hV{{\widehat V}}
\newcommand\hW{{\widehat W}}
\newcommand\hX{{\widehat X}}
\newcommand\hY{{\widehat Y}}
\newcommand\hZ{{\widehat Z}}

\newcommand\ha{{\widehat a}}
\newcommand\hb{{\widehat b}}
\newcommand\hc{{\widehat c}}
\newcommand\hd{{\widehat d}}
\newcommand\he{{\widehat e}}
\newcommand\hf{{\widehat f}}
\newcommand\hg{{\widehat g}}
\newcommand\hh{{\widehat h}}
\newcommand\hi{{\widehat i}}
\newcommand\hj{{\widehat j}}
\newcommand\hk{{\widehat k}}
\newcommand\hl{{\widehat l}}
\newcommand\hm{{\widehat m}}
\newcommand\hn{{\widehat n}}
\newcommand\ho{{\widehat o}}
\newcommand\hp{{\widehat p}}
\newcommand\hq{{\widehat q}}
\newcommand\hr{{\widehat r}}
\newcommand\hs{{\widehat s}}
\newcommand\htt{{\widehat t}}
\newcommand\hu{{\widehat u}}
\newcommand\hv{{\widehat v}}
\newcommand\hw{{\widehat w}}
\newcommand\hx{{\widehat x}}
\newcommand\hy{{\widehat y}}
\newcommand\hz{{\widehat z}}

\newcommand\tA{{\widetilde A}}
\newcommand\tB{{\widetilde B}}
\newcommand\tC{{\widetilde C}}
\newcommand\tD{{\widetilde D}}
\newcommand\tE{{\widetilde E}}
\newcommand\tF{{\widetilde F}}
\newcommand\tG{{\widetilde G}}
\newcommand\tH{{\widetilde H}}
\newcommand\tI{{\widetilde I}}
\newcommand\tJ{{\widetilde J}}
\newcommand\tK{{\widetilde K}}
\newcommand\tL{{\widetilde L}}
\newcommand\tM{{\widetilde M}}
\newcommand\tN{{\widetilde N}}
\newcommand\tOO{{\widetilde O}}
\newcommand\tP{{\widetilde P}}
\newcommand\tQ{{\widetilde Q}}
\newcommand\tR{{\widetilde R}}
\newcommand\tS{{\widetilde S}}
\newcommand\tTT{{\widetilde T}}
\newcommand\tU{{\widetilde U}}
\newcommand\tV{{\widetilde V}}
\newcommand\tW{{\widetilde W}}
\newcommand\tX{{\widetilde X}}
\newcommand\tY{{\widetilde Y}}
\newcommand\tZ{{\widetilde Z}}

\newcommand\ta{{\widetilde a}}
\newcommand\tb{{\widetilde b}}
\newcommand\tc{{\widetilde c}}
\newcommand\td{{\widetilde d}}
\newcommand\te{{\widetilde e}}
\newcommand\tf{{\widetilde f}}
\newcommand\tg{{\widetilde g}}
\newcommand\th{{\widetilde h}}
\newcommand\ti{{\widetilde i}}
\newcommand\tj{{\widetilde j}}
\newcommand\tk{{\widetilde k}}
\newcommand\tl{{\widetilde l}}
\newcommand\tm{{\widetilde m}}
\newcommand\tn{{\widetilde n}}
\newcommand\tio{{\widetilde o}}
\newcommand\tp{{\widetilde p}}
\newcommand\tq{{\widetilde q}}
\newcommand\tr{{\widetilde r}}
\newcommand\ts{{\widetilde s}}
\newcommand\tit{{\widetilde t}}
\newcommand\tu{{\widetilde u}}
\newcommand\tv{{\widetilde v}}
\newcommand\tw{{\widetilde w}}
\newcommand\tx{{\widetilde x}}
\newcommand\ty{{\widetilde y}}
\newcommand\tz{{\widetilde z}}

\newcommand\bA{{\overline A}}
\newcommand\bB{{\overline B}}
\newcommand\bC{{\overline C}}
\newcommand\bD{{\overline D}}
\newcommand\bE{{\overline E}}
\newcommand\bFF{{\overline F}}
\newcommand\bG{{\overline G}}
\newcommand\bH{{\overline H}}
\newcommand\bI{{\overline I}}
\newcommand\bJ{{\overline J}}
\newcommand\bK{{\overline K}}
\newcommand\bL{{\overline L}}
\newcommand\bM{{\overline M}}
\newcommand\bN{{\overline N}}
\newcommand\bO{{\overline O}}
\newcommand\bP{{\overline P}}
\newcommand\bQ{{\overline Q}}
\newcommand\bR{{\overline R}}
\newcommand\bS{{\overline S}}
\newcommand\bT{{\overline T}}
\newcommand\bU{{\overline U}}
\newcommand\bV{{\overline V}}
\newcommand\bW{{\overline W}}
\newcommand\bX{{\overline X}}
\newcommand\bY{{\overline Y}}
\newcommand\bZ{{\overline Z}}

\newcommand\ba{{\overline a}}
\newcommand\bb{{\overline b}}
\newcommand\bc{{\overline c}}
\newcommand\bd{{\overline d}}
\newcommand\be{{\overline e}}
\newcommand\bff{{\overline f}}
\newcommand\bg{{\overline g}}
\newcommand\bh{{\overline h}}
\newcommand\bi{{\overline i}}
\newcommand\bj{{\overline j}}
\newcommand\bk{{\overline k}}
\newcommand\bl{{\overline l}}
\newcommand\bm{{\overline m}}
\newcommand\bn{{\overline n}}
\newcommand\bo{{\overline o}}
\newcommand\bp{{\overline p}}
\newcommand\bq{{\overline q}}
\newcommand\br{{\overline r}}
\newcommand\bs{{\overline s}}
\newcommand\bt{{\overline t}}
\newcommand\bu{{\overline u}}
\newcommand\bv{{\overline v}}
\newcommand\bw{{\overline w}}
\newcommand\bx{{\overline x}}
\newcommand\by{{\overline y}}
\newcommand\bz{{\overline z}}

%
\newcommand{\AAA}{{\cal A}}
\newcommand{\BB}{{\cal B}}
\newcommand{\CC}{{\cal C}}
\newcommand{\DD}{{\cal D}}
\newcommand{\EE}{{\cal E}}
\newcommand{\FF}{{\cal F}}
\newcommand{\GG}{{\cal G}}
\newcommand{\HH}{{\cal H}}
\newcommand{\II}{{\cal I}}
\newcommand{\JJ}{{\cal J}}
\newcommand{\KK}{{\cal K}}
\newcommand{\LL}{{\cal L}}
\newcommand{\NN}{{\cal N}}
\newcommand{\MM}{{\cal M}}
\newcommand{\OO}{{\cal O}}
\newcommand{\PP}{{\cal P}}
\newcommand{\QQ}{{\cal Q}}
\newcommand{\RR}{{\cal R}}
\newcommand{\SS}{{\cal S}}
\newcommand{\TT}{{\cal T}}
\newcommand{\UU}{{\cal U}}
\newcommand{\VV}{{\cal V}}
\newcommand{\WW}{{\cal W}}
\newcommand{\XX}{{\cal X}}
\newcommand{\YY}{{\cal Y}}
\newcommand{\ZZ}{{\cal Z}}
\newcommand{\Id}{operatorname{Id}}

\newcommand{\tbullet}{$\bullet$}
\newcommand{\ot}{\leftarrow}
\newcommand{\carre}{\hfill$\Box$}
\newcommand{\carreb}{\hfill\rule{0.25cm}{0.25cm}}
%
%
\newcommand{\dontforget}[1]
{{\mbox{}\\\noindent\rule{1cm}{2mm}\hfill don't forget : #1
\hfill\rule{1cm}{2mm}}\typeout{---------- don't forget : #1 ------------}}
\newcommand{\note}[2]
{ \noindent{\sf #1 \hfill \today}

\noindent\mbox{}\hrulefill\mbox{}
\begin{quote}\begin{quote}\sf #2\end{quote}\end{quote}
\noindent\mbox{}\hrulefill\mbox{}
\vspace{1cm}
}
\newcommand{\rond}[1]     {\stackrel{\circ}{#1}}
\newcommand{\rondf}       {\stackrel{\circ}{\FF}}
\newcommand{\point}[1]     {\stackrel{\cdot}{#1}}

\newcommand\relatif{{\rm \rlap Z\kern 3pt Z}}

\title{\huge   Pseudoholomorphic discs near an elliptic point}
\author{Alexandre Sukhov* and Alexander Tumanov** }
\date{}
\maketitle

{\small 

{*}Universit\'e des Sciences et Technologies de Lille, Laboratoire
Paul Painl\'ev\'e,
U.F.R. de
Math\'ematique, 59655 Villeneuve d'Ascq, Cedex, France , 
 e-mail: sukhov@math.univ-lille1.fr

{**} Department of Mathematics, University of Illinois, 1409 West Green Street, Urbana, IL 61801,
USA, e-mail: tumanov@math.uiuc.edu}

\bigskip

Abstract. We prove the existence and  study the geometry of Bishop discs near an elliptic point 
of a real $n$-dimensional submanifold of an almost complex n-dimensional manifold.

MSC: 32H02, 53C15.

Key words: almost complex structure, generic manifold, elliptic point,
Bishop disc.

\bigskip

\section{Introduction}

E.Bishop \cite{Bi} proved that holomorphic discs with boundaries attached near an 
elliptic point  $p$ of a real $n$-dimensional submanifold $E$ of $\cc^n$ form a family smoothly 
depending on $(n-1)$ real parameters and that the boundaries of these discs fill a neighborhood 
of $p$ in $E$. A substantial progress in the study of the geometry of these discs (called Bishop discs) 
has been done by E.Bedford - B.Gaveau \cite{BeGa}, C.Kenig - S.Webster \cite{KeWe1,KeWe2} and other authors. On the other hand, M.Gromov 
in \cite{Gr} pointed out that Bishop discs still exist when $E$ is a submanifold of an almost
complex manifold; he used this result in his approach to the symplectic and contact geometry. In the case of almost complex dimension 2, a precise version of this result  
was established and substantially used by H.Hofer \cite{Ho}, Y.Eliashberg \cite{El}, R.Ye \cite{Ye} in various forms. The goal of the present paper is to give a 
complete proof of the existence of Bishop discs with boundaries attached near an elliptic point $p$ of a real $n$-dimensional submanifold of an almost 
complex manifold of complex dimension $n$. We prove that their boundaries fill a neighborhood of $p$ on $E$ similarly to the case of the standard complex structure. 
Our approach arises from the standard complex analysis techniques, in particular, we use  the methods of the  works  \cite{KS} and \cite{Tu}.

\section{Preliminaries}

\subsection{Almost complex manifolds.}
 Let $(M,J)$ be an almost complex manifold
with operator of complex structure $J$.
 Let $\D$ be the unit
disc in $\cc$ and $J_{st}$  the standard (operator of)  complex  structure
on $\cc^n$
for arbitrary  $n$. Let $f$ be
 a smooth map from $\D$ into $M$. We say that $f$ is {\it
 $J$-holomorphic}  if $df \circ J_{st} = J \circ df$. We call such a map $f$
a $J$-holomorphic disc and
 denote by   ${\mathcal O}_J(\D,M)$ the set
 of {\it $J$-holomorphic discs} in
 $M$.
We denote by ${\mathcal O}(\D)$ the space of usual holomorphic
 functions on $\D$.

The following lemma shows that an  almost complex manifold
$(M,J)$ can be  locally viewed  as the unit ball $\B$ in
$\cc^n$ equipped with a small almost complex
deformation of $J_{st}$. We shall repeatedly use this
observation in what follows.
\begin{e-lemme}
\label{lemma1}
Let $(M,J)$ be an almost complex manifold. Then for each $p \in
M$,  each  $\delta_0 > 0$, and  each   $k
\geq 0$
 there exist a neighborhood $U$ of $p$ and a
smooth coordinate chart  $z: U \longrightarrow \B$ such that
$z(p) = 0$, $dz(p) \circ J(p) \circ dz^{-1}(0) = J_{st}$,  and the
direct image $z_*(J) := dz \circ J \circ dz^{-1}$ satisfies
the inequality
$\vert\vert z_*(J) - J_{st}
\vert\vert_{\CC^k(\bar {\B})} \leq \delta_0$.
\end{e-lemme}
\proof There exists a diffeomorphism $z$ from a neighborhood $U'$ of
$p \in M$ onto $\B$ such that  $z(p) = 0$ and $dz(p) \circ J(p)
\circ dz^{-1}(0) = J_{st}$. For $\delta > 0$ consider the isotropic dilation
$d_{\delta}: t \mapsto \delta^{-1}t$ in $\cc^n$ and the composite
$z_{\delta} = d_{\delta} \circ z$. Then $\lim_{\delta \rightarrow
0} \vert\vert (z_{\delta})_{*}(J) - J_{st} \vert\vert_{\CC^k(\bar
{\B})} = 0$. Setting $U = z^{-1}_{\delta}(\B)$ for sufficiently small positive
$\delta$  we obtain the required result.

{\bf The operators $\partial_J$ and $\bar{\partial}_J$ .}

Let $(M,J)$ be an almost complex manifold. We denote by $TM$ the real
tangent bundle of $M$ and by $T_{\cc} M$ its complexification. Recall
that $T_{\cc} M = T^{(1,0)}M \oplus T^{(0,1)}M$ where
$T^{(1,0)}M:=\{ X \in T_{\cc} M : JX=iX\} = \{\zeta -iJ \zeta, \zeta \in
TM\},$
and $T^{(0,1)}M:=\{ X \in T_{\cc} M : JX=-iX\} = \{\zeta +
iJ \zeta, \zeta \in TM\}$.
 Let $T^*M$ be  the cotangent bundle of  $M$.
Identifying $\cc \otimes T^*M$ with
$T_{\cc}^*M:=Hom(T_{\cc} M,\cc)$ we define the set of complex
forms of type $(1,0)$ on $M$ as
$
T_{(1,0)}M=\{w \in T_{\cc}^* M : w(X) = 0, \forall X \in T^{(0,1)}M\}
$
and we denote  the set of complex forms of type $(0,1)$ on $M$ by
$
T_{(0,1)}M=\{w \in T_{\cc}^* M : w(X) = 0, \forall X \in T^{(1,0)}M\}
$.
Then $T_{\cc}^*M=T_{(1,0)}M \oplus T_{(0,1)}M$.
This allows us to define the operators $\partial_J$ and
$\bar{\partial}_J$ on the space of smooth functions  on
$M$~: for a  smooth complex function $u$ on $M$ we set $\partial_J u = du_{(1,0)} \in T_{(1,0)}M$ and $\bar{\partial}_Ju = du_{(0,1)}
\in T_{(0,1)}M$. As usual,
differential forms of any bidegree $(p,q)$ on $(M,J)$ are defined
by  exterior multiplication.

\subsection{Real submanifolds of an almost complex manifold} Let $E$ be a real submanifold of codimension $k$ in an almost complex manifold
$(M,J)$ of complex dimension $n$. Locally $E$ is defined as the zero-set of a smooth $\R^k$-valued function $r = (r_1,...,r_k)$ with 
$dr_1 \wedge ...\wedge dr_k \neq 0$. As usual, we denote by $H_pE$ the maximal complex subspace of the tangent space $T_pE$  that is 
 $H_pE = T_pE \cap J(p)T_pE$ and call it the holomorphic or the complex tangent space of $E$ at $p$. Recall that $E$ is called generic at $p$ 
if the complex dimension of $H_pE$ is equal to $n - k$. A real submanifold $E$ is called totally real at $p$ if $H_pE = \{ 0 \}$. In particular, 
any generic submanifold of real codimension $n$ is totally real.

 Let $E$ be a real $n$-dimensional submanifold in an almost complex manifold $(M,J)$.
Our considerations are local so we simply indentify $M$ with $\R^{2n}$ (or $\cc^n$ ) and view $J$ as a smooth matrix-valued 
function. In this paper we study the geometry of $E$ near a point $p \in E$ admitting a tangent vector $X \in T_pE$ sucht that 
$J(p)X \in T_pE$, that is a complex tangent direction that is $p$ is a singular point for the Cauchy-Riemann structure of $E$.
We consider generic singularities only that is the case where there exists precisely one tangent complex line in $T_pE$.
Our approach is based on a suitable choice of coordinates near $p$. First we can assume that $p = 0$ and $J(0) = J_{st}$.
Let vector $e_1= (1,0,...,0)$ be complex tangent to $E$ at the origin.  It follows by the classical Nijenhuis-Woolf theorem \cite{NiWo} 
that there exists a $J$-holomorphic disc $f:\D \longrightarrow \R^{2n}$ such that $f(0) = 0$ and $df(0)(\frac{\partial}{\partial (\Re \zeta)}) = e_1$.
After a local diffeomorphism with the identical linear part at the origin straightening $f$, we obtain that the map $\zeta \mapsto (\zeta,0,...,0)$
is $J$-holomorphic on $\D$. Denote the coordinates in $\cc^n$ by $Z = (w,z)$ with $z = (z_2,...,z_n)$. Then we may assume 
(possibly, after an additional change of coordinates preserving the $J$-holomorphicity of the coordinates line $w$ and a suitable choice of defining 
functions similarly to \cite{Bi}) that $E$ is defined near the origin by the equations 

\begin{eqnarray*}
& &\rho(Z) = z_2 - P(w) + R(Z) = 0,\\
& &r_3(Z) = \Re z_3 + h_3(Z) =0,\\
& &...\\
& &r_n(Z) = \Re z_n + h_n(Z) = 0
\end{eqnarray*}
where $R = O(\vert Z \vert^3 )$ and $h_j(Z) = O(\vert Z \vert^2)$. Here 
\begin{eqnarray*}
P(w) = w\overline w + \gamma(w^2 + \overline w^2)
\end{eqnarray*}
The point $p$ is called {\it elliptic} if in such a system of coordinates $\gamma \in [0,1/2[$. This notion is 
invariant. In what follows we denote by $r:= (\Re \rho, \Im \rho,r_3,...,r_n)$ the $\R^n$-valued local defining function of $E$.

\section{Bishop discs}

 Let $(M,J)$ be a
smooth almost complex
manifold of real dimension
$2n$ and $E$  a real submanifold of $M$ of real codimension
$m$. A $J$-holomorphic disc $f:\D \longrightarrow M$ continuous on
$\overline \D$ is called a {\it Bishop disc} if $f(b\D) \subset E$
(where $b\D$ denotes the boundary of $\D$).
 Our aim  is to prove the
existence and to describe certain classses of Bishop discs attached to $E$.
In this section we derive Riemann-Hilbert type boundary problems describing pseudoholomorphic Bishop discs near generic and elliptic points
of a real submanifold in an almost complex manifold. 

Since our considerations are local, we identify an almost complex manifold with $\cc^n$ and denote by $Z = (Z_1,...,Z_n)$ the standrard complex coordinates.
An almost complex structure $J$ may be viewed as a smooth real $(2n \times 2n)$-matrix valued function $J: Z \mapsto J(Z)$ on a neighborhood of 
the origin in $\cc^n$. We also may assume that $J(0) = J_{st}$. Then a smooth map 
$$Z:\zeta \mapsto Z(\zeta)$$
$$Z:\D \longrightarrow U$$
where $U$ is a small enough neighborhood of the origin in $\cc^n$ is $J$-holomorphic if and only if it satisfies the following PDE system
\begin{eqnarray*}
Z_{\overline\zeta} - A(Z) {\overline Z}_{\overline\zeta} = 0
\end{eqnarray*}
where $A(Z)$ is the
complex $n\times n$ matrix  of the operator whose composite  with complex
conjugation is equal to the endomorphism
$ (J_{st} + J_\delta(Z))^{-1}(J_{st} - J_\delta(Z))$
(which is an anti-linear operator with respect to the standard
structure $J_{st}$). By lemma \ref{lemma1} we can assume that the norm of $A$ is small enough.

\subsection{Bishop discs and Bishop's equation near a generic point} Denote by $Kf$ the Cauchy integral of a function $f$:

\begin{eqnarray*}
K(f)(\zeta) = \frac{1}{2\pi i}\int_{b \D}\frac{f(\tau) d\tau}{\tau -\zeta}
\end{eqnarray*}
for $\zeta \in \D$ and by $K_0f$ its principal value at boundary point:

\begin{eqnarray*}
K_0(f)(\zeta) = (v.p.) \frac{1}{2\pi i}\int_{b \D}\frac{f(\tau) d\tau}{\tau -\zeta}, \zeta \in \partial \D
\end{eqnarray*}

We will use also the following notation for the Cauchy-Green transform:

\begin{eqnarray*}
T(f)(\zeta) = \frac{1}{2\pi i}\int\int_{\D}\frac{f(\tau) d\tau \wedge d\overline \tau}{\tau -\zeta}
\end{eqnarray*}

Then for any  function $f$ of class $C^1$ on $\overline \D$ (the restriction $C^1$ can be replaced by some Sobolev's class, see \cite{Ve}) 
we have the generalized Cauchy formula:

\begin{eqnarray}
\label{Cauchy}
f = Kf + Tf_{\overline\zeta}
\end{eqnarray}

Denote also by $K_{+}f(\zeta)$ and $K_{-}f(\zeta)$ the limiting values of $Kf$ at $\zeta \in  b\D$ from inside and outside respectively.

Then we have the classical Plemejl-Sokhotski  formulae (at every point of $b\D$):
\begin{eqnarray}
\label{Plem1}
K_{\pm} f = K_0f \pm (f/2),
\end{eqnarray}
and 
\begin{eqnarray}
\label{Plem2}
K_+f - K_-f = f
\end{eqnarray}

In particular, (\ref{Cauchy}) and (\ref{Plem2}) imply that 

\begin{eqnarray}
\label{form1}
-T f(\zeta) = K_-f(\zeta), \zeta \in b \D
\end{eqnarray}

Recall also for any real function $u \in L^2(b\D)$ there exists a unique $v \in L^2(\bD)$ 
(the Hilbert transform of $u$) such that $u + iv$ is a boundary value of a holomorphic function $h$ from the Hardy space $H^2(b\D)$. 
The Hilbert transform
 may be written as a singular integral operator:

\begin{eqnarray*}
Hf(\zeta) = \frac{1}{2\pi i} (v.p) \int_{0}^{2\pi} \frac{e^{i\theta}  + \zeta}{e^{i\theta} - \zeta}f(e^{i\theta})d\theta
\end{eqnarray*}
The classical fact is that for a non-integer $k> 0$, $H$ is a continuous linear operator on $C^k(b\D)$. 
Denote also by $P_0f$ the value of the Poisson integral of $f$ at the origin:

\begin{eqnarray*}
P_0f = \frac{1}{2\pi}\int_0^{2\pi}f(e^{i\theta})d\theta
\end{eqnarray*}

Integrating the identity $2\tau(\tau - \zeta)^{-1} = (\tau + \zeta)(\tau - \zeta)^{-1} + 1$ we obtain that

\begin{eqnarray*}
\label{form2}
K_0f = (i/2)Hf + (1/2)P_0f
\end{eqnarray*}
on the boundary of $\D$. Now (\ref{Plem1}) gives

\begin{eqnarray}
\label{form3}
f(\zeta) = iHf(\zeta) + P_0f - 2K_-f(\zeta), \zeta \in b\D
\end{eqnarray}

Suppose that $f$ is of class $C^{\alpha}(\overline \D)$ with $\alpha > 0$ non-intergral. Let also  $f^* = u^* + iv^*$
be a function on $b\D$ of class $C^{\alpha}$.

\begin{e-lemme}
Assume that $f$ satisfies the following system of singular integral equations:
\begin{eqnarray}
\label{form4}
f = Kf^* + T f_{\overline\zeta}, \,\,\, \zeta \in \D
\end{eqnarray}
\begin{eqnarray}
\label{form5}
v^* = Hu^* + P_0v^* + 2\Im T f_{\overline\zeta}, \,\,\, \zeta \in b\D
\end{eqnarray}
Then $f^* = f\vert_{ b\D}$.
\end{e-lemme}
\proof (\ref{form4}) implies that $K_+f^* = K_+f$. In view of (\ref{form1}) we obtain $K_+f = Kf^*$.
Furthermore, (\ref{form5}) implies $-2\Im T f_{\overline\zeta} = \Im (-v^* + Hu^* + P_0v^*) = \Im (iHf^* + P_0f^* - f^*) = 2\Im K_-f^*$.
Therefore, $K_f = K_-f^*$ (recall that the value of the Cauchy type integral at the infinity is equal to $0$). Thus, 
$K_+f = K_+f^*$ and $K_-f = K_-f^*$ which implies the statement of lemma.

Now we may write the Bishop equation. We assume that a real submanifold $E$ is given by the following equations near the origin:

\begin{eqnarray*}
y = h(x,w)
\end{eqnarray*}
with 
\begin{eqnarray*}
h(0,0) = 0, dh(0,0) = 0,
\end{eqnarray*}
where we use the notation $Z = (z,w) \in \cc^{n-m} \times \cc^m$, $z = x + iy$, $w = u + iv$.  
We also write a map $Z$ in the form $Z(\zeta) = (z(\zeta),w(\zeta))$. We will represent the matrix $A$ in the  form 

$$
A= \left(
\begin{array}{cll}
P & & R\\
L & & N 
\end{array}
\right).
$$
with $(m\times m)$-bloc $P$.

Then a smooth disc $Z$ is $J$-holomorphic 
and its boundary is attached to $E$ if and only if it satisfies the following system of singular integral equations:

\begin{displaymath}
(RH_0): \left\{ \begin{array}{lll}
z = Kz^* - T(P(Z){\overline z}_{\overline\zeta} + R(Z){\overline w}_{\overline\zeta}),\,\,\,\zeta \in \D\\
y^* = H h(y^*,w^*) + y_0 - 2\Im T (P(Z) {\overline z}_{\overline\zeta} + R(Z) {\overline w}_{\overline\zeta}),\,\,\, \zeta \in b\D\\
w = Kw^* - T (L(Z){\overline z}_{\overline\zeta} + N(Z){\overline w}_{\overline\zeta}),\,\,\, \zeta \in \D\\
v^* = H u^* + v_0 - 2\Im T (L(Z) {\overline z}_{\overline\zeta} + N(Z) {\overline w}_{\overline\zeta}), \,\,\, \zeta \in b\D
\end{array} \right \}
\end{displaymath}
where $Z^* = Z \vert_{b\D}$.
Thus, this non-linear boundary Riemann-Hilbert type problem may be viewed as an almost complex analog of the standard Bishop equation.   
Applying the implicit function theorem we easily obtain that given non-integral $\alpha >0$, real constants $y_0$, $v_0$ and a real function $u^* \in C^\alpha(b\D)$ 
there exists a unique solution of the system $(RH_0)$.

 We point out that the existence of Bishop discs attached to a generic submanifold of an almost complex manifold was proved in \cite{KS} 
by a different method. The present method is substantially more explicit and allows to write Bishop's equation as a singular integral equation 
similarly to the integrable case. This leads to a more convenient and explicit parametrization of the Banach space of Bishop's discs and can 
be useful for other applications. In particular, it follows immediately from the equation $(RH_0)$ that Bishop's discs depend smoothly on 
deformations of the almost complex sctructure $J$.

\subsection{Bishop discs and Bishop's equation near an ellpitic point} Consider a real $n$-dimensional submanifold $E$ in an almost complex $n$-dimensional 
manifold $(M,J)$ and let $p$ be an elliptic point of $E$. 
We assume without loss of generality that $E$ is a submanifold in $\cc^n$ and the origin is an elliptic point of $E$; moreover, $J(0) = J_{st}$.
Furthemore, choosing local coordinates similarly to the previous section, we assume that 
the map $\zeta \mapsto (\zeta,0,...,0)$
is $J$-holomorphic on $\D$. Denote as above the coordinates in $\cc^n$ by $Z = (w,z)$ with $z = (z_2,...,z_n)$. Then we may assume 
 that $E$ is defined near the origin by the equations 

\begin{eqnarray*}
& &\rho(Z) = z_2 - P(w) + R(Z) = 0,\\
& &r_3(Z) = \Re z_3 + h_3(Z) =0,\\
& &...\\
& &r_n(Z) = \Re z_n + h_n(Z) = 0
\end{eqnarray*}
where $R = O(\vert Z \vert^3 )$ and $h_j(Z) = O(\vert Z \vert^2)$. Here 
\begin{eqnarray*}
P(w) = w\overline w + \gamma(w^2 + \overline w^2)
\end{eqnarray*}
and $\gamma \in [0,1/2[$. As above,  we denote by $r:= (\Re \rho, \Im \rho,r_3,...,r_n)$ the $\R^n$-valued local defining function of $E$.

Then a smooth map $f$
defined on $\D$ and continuous on
$\overline\D$ is a Bishop disc if and only if it satisfies the
following non-linear boundary problem of the
Riemann-Hilbert type for the
quasi-linear operator $\overline\partial_J$:

\begin{displaymath}
(RH): \left\{ \begin{array}{ll}
\overline\partial_J f(\zeta) = 0, \zeta \in \D\\
r(f)(\zeta) = 0, \zeta \in b\D
\end{array} \right.
\end{displaymath}

Consider the non-isotropic dilations
$\Lambda_{\delta}: Z = (w,z) \mapsto Z' = (\delta^{-1/2}w, \delta^{-1}z)$. In the new
$Z$-variables  (we drop the primes) the image
$E_{\delta} = \Lambda_{\delta}(E)$ is
defined by the equation $r_\delta (Z):=\delta^{-1}r((\Lambda_\delta)^{-1} Z) = 0$. 
Then the vector-function $r_\delta$ tends to the function 
$$r^0(Z) = (\Re z_2 - \Re P(w), \Im z_2 - \Im P(w), \Re z_3 + Q_3(w),...,\Re z_n + Q_n(w))$$ 
as $\delta \longrightarrow 0$. Here $Q_j$ are some real homogeneous of degree 2 polynomials in $w$.
Hence 
the manifolds $E_\delta$ approach the model quadric manifold $E_0=\{
r^0(Z) = 0\}$ which explicitly is given by the equations

\begin{eqnarray*}
& &\rho^0(Z) = z_2 - P(w)  = 0,\\
& &r_3^0(Z) = \Re z_3 - Q_3(w) =0,\\
& &...\\
& &r_n^0(Z) = \Re z_n - Q_n(w) = 0
\end{eqnarray*}

Consider the transported structures $J_\delta := (\Lambda_\delta)_*(J) = d\Lambda_\delta \circ J \circ (d\Lambda_\delta)^{-1}$.
The following statement is similar to \cite{KS}.

\begin{e-lemme}
For any positive $k$ and any compact subset $K \subset \cc^n$ we have $\parallel J_\delta - J_{st} \parallel_{C^k(K)} \longrightarrow 0$ 
as $\delta \longrightarrow 0$.
\end{e-lemme}

\proof Consider the Taylor expansion of $J(Z)$ near the origin: $J(Z)
= J_{st} + L(Z) + R(Z)$ where $L(Z)$ is  the linear part of the expansion
 and $R(Z)= O(\vert Z \vert^2)$. Clearly, $\Lambda_{\delta} \circ
R(\Lambda_{\delta}^{-1}(Z)) \circ \Lambda_{\delta}^{-1}$ converges to $0$
as $\delta \longrightarrow 0$. Denote by $L_{kj}^{\delta}(Z)$
(respectively, by $L_{kj}(Z)$) an entry
of the real  matrix $\Lambda_{\delta} \circ
L(\Lambda_{\delta}^{-1}(Z)) \circ \Lambda_{\delta}^{-1}$
(respectively, of $L(Z)$). For $k = 3,...,2n$, $j = 1,2$
we have $J_{kj}^{\delta}(w,z) = \delta^{-1/2}L_{kj}(\delta^{1/2}
w,\delta z) \longrightarrow L_{kj}(w,0)$ as $\delta \longrightarrow 0$. However, in the
 coordinate system fixed above the map $\zeta \longrightarrow (\zeta,0,..,0)$ is $J$-holomorphic,
so 
$L_{kj}(w,0) = 0$ for $k,j =1,2$. This implies that 
$L_{kj}^{\delta}$ approaches  $0$ for all $k,j$ (for other values of $k,j$ this is obvious). This gives us the
result of the lemma.

Consider a $J_\delta$-holomorphic disc
$Z: \D \longrightarrow (\B^n,J_\delta)$. The
$J_\delta$-holomorphy condition $J_\delta(Z) \circ dZ = dZ \circ J_{st}$ can be
written in the following form:

\begin{eqnarray}
\label{CR}
Z_{ \overline\zeta} - A_{J,\delta}(Z)
{ \overline Z}_{ \overline \zeta}  
= 0
\end{eqnarray}
where as above $A_{J,\delta}(Z)$ is the
complex $n\times n$ matrix  of the operator whose composite  with complex
conjugation is equal to the endomorphism
$ (J_{st} + J_\delta(Z))^{-1}(J_{st} - J_\delta(Z))$
(which is an anti-linear operator with respect to the standard
structure $J_{st}$). Hence the entries of the matrix $A_{J,\delta}(Z)$
are smooth functions of $\delta,Z$ vanishing identically in $Z$
for $\delta = 0$.

Using the Cauchy-Green transform $T_{CG}$ defined above we may write   $\overline\partial_J$-equation (\ref{CR}) as follows:

$$ [ Z - T_{CG} A_{J,\delta}(Z){ \overline Z}_{ \overline \zeta} ]_{\overline\zeta} = 0$$

 This is well known \cite{Ve} that  the Cauchy-Green
 transform is a continuous linear operator from $C^k(\overline \D)$ into
 $C^{k+1}(\overline \D)$ (recall that $k$ is nonintegral). Hence  the operator

$$\Phi_{J,\delta}: Z \longrightarrow W =  Z - T_{CG}A_{J,\delta}(Z){ \overline Z}_{ \overline \zeta}$$
takes the space   $C^{k}(\overline\D)$  into itself.
Thus, $Z$ is $J_\delta$-holomorphic if
and only if $\Phi_{J,\delta}(Z)$ is holomorphic (in the usual sense)
on $\D$.
For sufficiently small positive $\delta$  this is an
invertible operator on a neighborhood  of zero in $C^k(\overline \D)$
which  establishes a one-to-one correspondence between the
sets of $J_\delta$-holomorphic and holomorphic discs  in $\B^n$.

These considerations allow us to replace the non-linear Riemann-Hilbert problem
(RH) by  Bishop's equation

\begin{eqnarray}
\label{Bishop2}
r_\delta(\Phi_{J,\delta}^{-1}(W))(\zeta) = 0, \zeta \in b\D
\end{eqnarray}
for  an unknown  holomorphic function $W$ in the unit  disc.

If $W$ is a solution of the boundary  problem (\ref{Bishop2}), then $Z = \Phi_{J,\delta}^{-1}(W)$ is a Bishop disc with  boundary
attached to $E_\delta$. Since the manifold $E_\delta$ is
biholomorhic via  non-isotropic dilations to the initial manifold $E$,
the solutions of the equation (\ref{Bishop2}) allow to describe
the Bishop discs attached to $E$.

\subsection{Analysis of the Bishop equation and the geometry of Bishop discs}

We begin with the description of Bishop discs attached to the model manifold $E_0$ in $\cc^n$ with the standrard structure $J_{st}$.
They are the solutions of the boundary problem (\ref{Bishop2}) for $\delta = 0$.

For $r > 0$  consider the ellipse $D_r: = \{ w \in \cc: P(w) < r \}$  and denote by $w_r$ the biholomorphism 
$w_r = w(r,\bullet):\D \longrightarrow D_r$ satisfying $w_r(0) = 0$. Then $P \circ w_r \vert_{ b\D} \equiv  r$ and we set $z_2(r,\zeta)
\equiv r$. This $J_{st}$ holomorphic disc admits the lift to a Bishop disc if we set 

\begin{eqnarray*}
z_j(r,\zeta) = S Q_j(\zeta) + ic_n, j=3,...,n
\end{eqnarray*}
where $c_3,...,c_n$ are real constans and $S$ denotes  the Schwarz integral:
\begin{eqnarray*}
S f(\zeta) =  \frac{1}{2\pi i}\int_{\vert \tau  \vert = 1} \frac{\tau + \zeta}{\tau - \zeta}f(\tau) \frac{d\tau}{\tau}
\end{eqnarray*}

We obtain a family of $J_{st}$-holomorphic discs $(w(r,c)(\bullet),z(r,c)(\bullet))$ smoothly depending on $(n-1)$ real parameters $(r,c_3,...,c_n)$, $r > 0$.
Their boundaries are disjoint  and form a foliation of a neighborhood of the real $(n-2)$ dimensional submanifold  
$$\Sigma:= \{ w = 0, z_2 = 0, \Re z_3 =... = \Re z_n = 0 \}$$
  in $E$. For $r = 0$
these discs degenerate to the constant map $\zeta \mapsto (0,0,ic_3,...,ic_n)$. Their images  form a real  $(n+1)$ submanifold $\tilde E$  
such that $(\tilde E,E)$ is a smooth manifold with boundary outside the points of $\Sigma$. 

We claim that in the general case of an almost complex structure $J$ the Bishop discs have similar properties. Indeed,   Let $w=0$, $z_2 = \delta$, $z_j = \delta$,
$j=3,...,n$ for $\delta > 0$ be a point on the real ``normal'' to $E$. After the non-isotropic dilation $\Lambda_\delta$ the image of this point coincides with  
$p^0:= (0,1,1,...,1)$. There exists a unique Bishop disc $Z^0$ of the described above family centered at $p^0$; it corresponds to the parameters $r = 1$ and $c = 0$. 
The parametrizing map $F:(r,c,\zeta) \longrightarrow (w(r,c)(\zeta),z(r,c)(\zeta))$ has the maximal rank when the parameters $(r,c)$ 
are in a neighborhood $U$  of the point $(1,0)$ and $\zeta \in \D$. Applying the implicit function theorem to the operator equation (\ref{Bishop2}) for $\delta > 0$ 
small enough we obtain 
a holomorphic solution $W(r,c)(\bullet)$ smoothly depending on $(n-1)$ real parameters $(r,c)$. Then $Z(r,c) = \Phi_{J,\delta}^{-1}(W)$ is a family of $J$-holomorphic Bishop discs with boundaries attached to $E_\delta$. This family is a small deformation of the above Bishop disc for $E_0$. Since $\delta$ is small, the parametrizing map 
$(r,c,\zeta) \mapsto Z(r,c)(\zeta)$ has the maximal rank. So these discs form a real smooth $(n+1)$ manifold $\tilde E_\delta$ with smooth boundary. Furthermore, 
this manifold is foliated by the Bishop discs $Z(r,c)(\bullet)$  and a neigborhood of $Z^0(b\D)$ in $E$ is an open subset of the  boundary of $\tilde E_\delta$.

We proved the following statement.

\begin{e-theo}
Let $p$ be an elliptic point of a real $n$ dimensional submanifold $E$ of an almost complex manifold $(M,J)$. Given positive $k$ there exists a family 
of  $J$-holomorphic Bishop discs for $E$ $C^k$-smoothly depending on $(n-1)$ real parameters. Thiese discs foliate a real $(n+1)$-dimensional submanifold 
$\tilde E$ such that $(\tilde E, E)$ is a $C^k$ smooth manifold with boundary outside a real $(n-2)$-dimensional submanifold $\Sigma$ in $E$.
\end{e-theo} 

As we pointed out, in the case of dimension 2 a similar statement is due to H.Hofer \cite{Ho}, Y.Eliashberg \cite{El} and  R.Ye \cite{Ye}.
It seems natural to study the questrions related to the regularity of the manifold $\tilde E$ near the point of $\Sigma$ similarly to C.Kenig - S.Webster results \cite{KeWe1,
KeWe2} in the case of an integrable complex structure.

{\it Aknowledgements.} This paper was done while the first named author visited the University of Illinois at Urbana-Champaign during the Spring semester of 2005.
He thanks this institution for excellent conditions for work.

\end{document}